\documentclass[article,12pt]{elsarticle}
\makeatletter
\def\ps@pprintTitle{%
 \let\@oddhead\@empty
 \let\@evenhead\@empty
 \def\@oddfoot{\centerline{\thepage}}%
 \let\@evenfoot\@oddfoot}
\makeatother
\newtheorem{property}{Property}
\newtheorem{proposition}{Proposition}
\newtheorem{theorem}{Theorem}
\newtheorem{lemma}{Lemma}

\usepackage{mathtools, amsmath, amssymb, graphicx, amsfonts}
\usepackage{xcolor} %to color comments
\usepackage{bbm} %to use indicator function
\usepackage{float} %to use [h!] with figures
\usepackage{amsopn} %standard?
\usepackage{mathrsfs} %for \mathscr
\usepackage{comment} %to comment out sections
\usepackage[normalem]{ulem} %to strikout text

\usepackage{pst-node}
\usepackage{tikz-cd} 

\usepackage{wrapfig} %to wrap around figure
\usepackage{lineno} %for line numbers

\newcommand{\tr}{\operatorname{tr}}

\journal{Linear Algebra and its Applications}

\begin{document}

\begin{frontmatter}

\title{Effective resistance is more than distance: \\Laplacians, Simplices and the Schur complement}
\tnotetext[t1]{The author was supported by The Alan Turing Institute under the EPSRC grant EP/N510129/1.}
\author{Karel Devriendt \\\emph{University of Oxford, Oxford, UK\\Alan Turing Institute, London, UK}}

%%%%%%%%%%%%%%
%% ABSTRACT %%
%%%%%%%%%%%%%%
\begin{abstract}
This article reviews and discusses a geometric perspective on the well-known fact in graph theory that the effective resistance is a metric on the nodes of a graph. The classical proofs of this fact make use of ideas from electrical circuits or random walks; here we describe an alternative approach which combines geometric (using simplices) and algebraic (using the Schur complement) ideas. These perspectives are unified in a matrix identity of Miroslav Fiedler, which beautifully summarizes a number of related ideas at the intersection of graphs, Laplacian matrices and simplices, with the metric property of the effective resistance as a prominent consequence.
\end{abstract}

\end{frontmatter}

%\linenumbers

%%%%%%
%BODY
%%%%%%
\section{Introduction}\label{Section: introduction}
The Laplacian matrix was first formulated (implicitly) by Gustav Kirchhoff in the context of electrical circuits in \cite{Kirchhoff, Kirchhoff_translated}, where it captures the linear relation between voltages and currents in a circuit of resistors. Results such as Kirchhoff's Matrix Tree Theorem however, which states that the number of spanning trees of a given connected graph can be found as the product of the non-zero Laplacian eigenvalues divided by the number of nodes, highlight that the Laplacian is a fundamental object in the study of graphs independent of the context of electrical circuits. The same story holds for the effective resistance, which made its way from a concept in electrical circuit theory to an important graph property after discoveries such as its relation to random walks \cite{Nash-Williams, Doyle, Palacios_resistances_and_randomwalks}, its role in the famous problem of ``dissecting the rectangle into squares" \cite{Tutte}, its function as a graph metric \cite{Klein, Qiu_Clustering_commute_times} and a robustness measure \cite{Ghosh, Ellens, Tyloo}, and more recently its role in ``spectral sparsification" of graphs \cite{Spielman}. The effective resistance also plays an important role in chemical graph theory, where it is used in formulating the so-called Kirchhoff index \cite{Klein, Palacios_Kirchhoff, Gutman, Bendito_Carmona}. These parallel histories are no coincidence, but reflect the intimate connection between Laplacians and effective resistances -- a connection perhaps best captured by a beautiful result due to Miroslav Fiedler that represents their relation in a single matrix identity \cite[Thm. 1.4.1]{Fiedler_book} (see also Section \ref{section: effective resistance and fiedler}).
\\
As mentioned, one of the important properties of the effective resistance is that it determines a metric between the nodes of a graph. Together with the geodesic distance, this resistance metric is probably one of the most natural notions for distances on a graph with the additional advantage of an exact algebraic expression based on the graph structure, and which can be calculated $\epsilon$-close in $O(m/\epsilon^2)$ time for $m$-link graphs \cite{Spielman}. The metric property of the effective resistance was discovered independently by Gvishiani and Gurvich in \cite{Gvishiani_Gurvich} and Klein and Randi\'{c} in \cite{Klein} using simple arguments from electrical circuit theory, and alternative proofs follow from the equivalence between effective resistances and commute times of random walks \cite{Gobel_Jagers_commute, Chandra, Nash-Williams}. While very concise, these proof strategies leave several key features of the effective resistance obscured. 
\\
Firstly, the usual discussions about the effective resistance do not make use of the fact that the square root of the effective resistance is a \emph{Euclidean metric} which associates to each weighted graph a simplex, and conversely. This fact, discovered independently by Fiedler in the context of simplex geometry \cite{Fiedler_aggregation} and by Sharpe and Moore in the context of electrical circuit theory \cite{Sharpe}, provides a distinct \emph{geometric} perspective for the study of effective resistances and, by extension, graphs. A second important result about effective resistances is the fact that the Schur complement (a certain operation that reduces a matrix to a smaller matrix on a subset of its columns and rows) maps Laplacian matrices to Laplacian matrices \emph{and} corresponds to a map from simplices to simplices. This result gives a complementary \emph{algebraic} perspective on the effective resistance. Finally, Fiedler's identity between Laplacian and effective resistance matrices unifies the algebraic and geometric viewpoints into a single concise matrix identity. While this result reflects several key properties of the effective resistance, it does not seem to be widely known and is rarely mentioned in the context of effective resistances or Laplacian matrices. 
\\
In this article, we present a self-contained derivation of the geometric (related to simplices) and algebraic (related to the Schur complement) characteristics of the effective resistance. As a particular application, we show how this setup admits an elegant proof for the triangle inequality of the effective resistance; however, as the title suggests, by the time we arrive at this metricity result it will be clear that this is indeed just one of the many qualities of the effective resistance. While most of the presented results have been described before by Fiedler \cite{Fiedler_book} or in the context of distance geometry, we believe that a unified and concentrated exposition of these ideas might be necessary for a wider understanding of the results and to promote their adoption by other researchers. In our conclusion, we furthermore suggest one possible path for future research that would consist of `categorifying' the ideas presented in this article and continuing further investigations in the more abstract realm of category theory.
\\
Our approach necessarily leaves many details unexplored and for additional results we refer the readers to \cite{Fiedler_book,krl_simplex} for the graph-simplex correspondence, \cite{Haynesworth, Dorfler_kron} for the Schur complement and \cite{Doyle, Klein, Ellens, Gutman, Bendito_Carmona, Gvishiani_Gurvich, Slepian_book, Biggs_potential, Bapat_book, Dorfler_electrical,  Sun} for electrical circuit theory and the effective resistance.
\\~
\\
The remainder of this article is organized as follows: in Section \ref{section: Graphs, Laplacians and Simplices} we introduce Laplacian matrices, simplices and effective resistances and discuss how these different objects are related; the main relations are summarized in Theorem \ref{Theorem: bijection between Laplacians and Simplices} and Theorem \ref{Theorem: Fiedlers identity}. Section \ref{Section: faces and Schur} then discusses how different objects of the same type are related, i.e. how faces of simplices are again simplices and how, correspondingly, the Schur complement maps Laplacian matrices to Laplacian matrices. In Section \ref{Section: geometric proof of resistance distance} finally, we present a simple geometric proof of the distance property of the effective resistance, Theorem \ref{Theorem: resistance is distance}, highlighting the utility and value of the earlier developed results. As an outlook on future work, we conclude the article with a brief description of a categorical perspective on the results we have introduced.
%%%
%
%%%
\\~\\
\section{Graphs, Laplacians and Simplices}\label{section: Graphs, Laplacians and Simplices}
\subsection{Laplacian matrices}\label{section: Laplacian matrices}
\begin{wrapfigure}{r}{0.4\textwidth}
    \fbox{
    \includegraphics[width=0.34\textwidth]{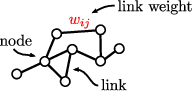}}%was scale=0.8
\end{wrapfigure}
\noindent We start by defining the basic objects of interest. A \emph{weighted graph} $G=(\mathcal{N},\mathcal{L},w)$ consists of a set of $n$ nodes $\mathcal{N}$ and a set of links $\mathcal{L}$ which connect (unordered, distinct) pairs of nodes\footnote{Nodes and links are often called vertices and edges in the graph theory literature. Here, these terms are reserved for the vertices and edges of a simplex.}, and positive weights $w$ defined on the links; we write $(i,j)\in\mathcal{L}$ for a link\footnote{An undirected link is sometimes denoted as $\lbrace i,j\rbrace$ instead of a tuple, to distinguish it from a directed link.} between $i$ and $j$, and $w_{ij}$ for its weight.
We assume the graph to be finite\footnote{We also assume $n\geq 2$ to avoid the exceptional case of the trivial graph.} ($n<\infty$) and connected, i.e. with a path between any pair of nodes. It is often more practical to represent the graph structure as a matrix; here, we work with the \emph{Laplacian matrix} of a graph $G$ which is the $n\times n$ matrix $Q$ with entries \cite{Merris, Mohar}
\begin{equation}\label{eq: definition Laplacian matrix}
(Q)_{ij} = \begin{cases}
d_i&\text{~if~} j=i\\
-w_{ij}&\text{~if~} (i,j)\in\mathcal{L}\\
0&\text{~otherwise,}
\end{cases}
\end{equation}
where $d_i$ is the \emph{degree} of a node $i$, equal to the total weight of all links containing $i$ as $d_i=\sum_{(i,j)\in\mathcal{L}}w_{ij}$. The properties of a (finite, connected and positively) weighted graph $G$ translate to properties of the Laplacian $Q$ as follows; the Laplacian matrix is/has:
\begin{align*}
    &\text{(i) symmetric}\\
    &\text{(ii) finite non-positive off-diagonal entries}\\
    &\text{(iii) zero row and column sums}\\
    &\text{(iv) irreducible}
\end{align*}
where irreducibility means that the matrix can not be block diagonalized by any permutation of the rows and columns. If we take properties (i)-(iv) as the definition of a Laplacian matrix, then expression \eqref{eq: definition Laplacian matrix} determines a bijection between Laplacian matrices and weighted graphs. We will continue with the Laplacian description, keeping in mind that any Laplacian corresponds to a weighted graph where the set of row/column indices $\mathcal{N}$ corresponds to the nodes of the graph and the set of non-zero off-diagonal entries $\mathcal{L}$ to links of the graph. 
\\
From the Laplacian properties above, a number of spectral properties of the Laplacian follow quite straightforwardly (see Proof of Proposition \ref{proposition: equivalent characterizations for Laplacian matrices} below); the Laplacian matrix is:
\begin{align*}
    &\text{(i)}_\sigma \text{~positive semidefinite}\\
    &\text{(ii)}_\sigma \text{~has a single zero eigenvalue,}\\
    &\text{(iii)}_\sigma \text{~with corresponding constant eigenvector}
\end{align*}
We will write the constant all-one vector as $u=(1,\dots,1)^T$. These spectral ($\sigma$ for spectral) properties\footnote{We remark that the positive semidefinite condition (i)$_{\sigma}$ implies both non-negative eigenvalues as well as symmetry of the matrix, even though the latter is not a spectral property.} suggest an interesting alternative definition of the Laplacian, with complementary information on the structure of Laplacian matrices:
\begin{proposition}\label{proposition: equivalent characterizations for Laplacian matrices}
The following characterizations for a matrix $Q$ are equivalent:
\begin{enumerate}\itemsep-1em 
    \item $Q$ is a Laplacian matrix\\
    \item $Q$ satisfies properties (i)-(iv)\\
    \item $Q$ satisfies properties (i)$_\sigma$-(iii)$_\sigma$ and  (ii)
\end{enumerate}
\end{proposition}
\textbf{Proof:} $1.\Leftrightarrow 2.$ follows from definition \eqref{eq: definition Laplacian matrix} of the Laplacian matrix and the fact that a connected graph corresponds to an irreducible Laplacian matrix. $1. \Rightarrow 3.$ From its definition, we have that a quadratic product with the Laplacian can be written as $x^TQx=\sum_{(i,j)\in\mathcal{L}}w_{ij}(x_i-x_j)^2\geq 0$ for any vector $x\in\mathbb{R}^n$, and thus all eigenvalues must be non-negative. Moreover, equality (i.e. zero eigenvalue) only holds when $x_i=x_j$ for all linked nodes $(i,j)\in\mathcal{L}$ and thus, by connectedness of the corresponding graph, for all nodes $i,j\in\mathcal{N}$ (i.e. constant eigenvector corresponding to the zero eigenvalue). $3.\Rightarrow 2.$ As $Q$ is positive semidefinite, it must be symmetric so that (i) is satisfied. Since the constant eigenvector has corresponding eigenvalue zero, we have that $\sum_j (Q)_{ij}=0$ for all $i$ so that (iii) is satisfied. Furthermore, if we assume $Q$ were reducible and could thus be written in the form $Q=\left(\begin{smallmatrix}Q_1&0\\0&Q_2 \end{smallmatrix}\right)$ then a vector $u'=(0,u)$ would have $Qu'=0$. However, since $Q$ is positive semidefinite with a single constant zero eigenvector $u\neq u'$ this is not possible, hence $Q$ is irreducible and (iv) holds. \hfill$\square$
\\
From Proposition \ref{proposition: equivalent characterizations for Laplacian matrices} it follows that the Laplacian has a spectral decomposition of the form
\begin{equation}\label{eq: Laplacian eigendecomposition}
Q=\sum_{k=1}^{n-1}\mu_kz_kz_k^T
\end{equation}
with real eigenvalues $\mu_k>0$ and normalized eigenvectors $z_k$ that satisfy the eigenequation $Qz_k=\mu_kz_k$, and where the zero eigenvalue $\mu_n=0$ and corresponding constant eigenvector $z_n=u/\sqrt{n}$ are omitted. The eigenvectors $\lbrace z_k\rbrace_{k=1}^n$ form an orthonormal basis for $\mathbb{R}^n$. However, as illustrated by the third characterization in Proposition \ref{proposition: equivalent characterizations for Laplacian matrices}, this eigendecomposition is not sufficient for $Q$ to be a Laplacian matrix as it does not constrain the off-diagonal entries to be non-positive; there is no simple `spectral fingerprint' that guarantees this sign property. This particular nature of the off-diagonal sign constraints will be discussed more later. 
\\
Another consequence of Proposition \ref{proposition: equivalent characterizations for Laplacian matrices} and decomposition \eqref{eq: Laplacian eigendecomposition} is that we can define the Moore-Penrose \emph{pseudoinverse Laplacian} $Q^\dagger$ as the inverse of $Q$ in the space orthogonal to the constant vector $u$, see for instance \cite{Fouss_Saerens_book}. In other words, such that $Q^\dagger Q = QQ^\dagger=I-uu^T/n$ which is a projector on the subspace orthogonal to $u$. More precisely, we can define the pseudoinverse Laplacian\footnote{Other notions of matrix pseudoinversion exist, but here we use `pseudoinverse' and the $\dagger$-superscript to refer to the Moore-Penrose pseudoinverse.} via its spectral decomposition as
$$
Q^\dagger = \sum_{k=1}^{n-1}\mu_k^{-1}z_kz_k^T
$$
which shows that $Q^\dagger$ is also positive semidefinite (i)$_\sigma$ with a single zero eigenvalue (ii)$_\sigma$ and constant eigenvector (iii)$_\sigma$. This spectral decomposition shows that properties (i)$_{\sigma}$-(iii)$_{\sigma}$ are always conserved under taking the Moore-Penrose pseudoinverse of the Laplacian matrix.
\\
\textit{Remark:} When a graph is not connected but consists of $\beta_0$ components, the corresponding Laplacian matrix will have a $\beta_0$-dimensional zero eigenspace spanned by eigenvectors which are piecewise constant on the components. In the language of algebraic topology, this corresponds to the fact that the zeroth Betti number (i.e. number of components) equals the dimension of the zeroth homology group, which in turn equals the dimension of the kernel of the Hodge Laplacian (i.e. our Laplacian), if our (weighted) graph is interpreted as a simplicial $1$-complex \cite{Lek-Heng_Lim_Hodge_Laplacians, Biggs_potential, Hansen_Spectral_sheaf_theory}.
%%%
%
%%%
\subsection{Simplices}\label{section: simplices}
\begin{figure}[h!]
    \centering\fbox{
    \includegraphics[width=0.9\textwidth]{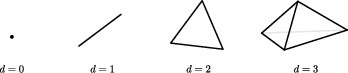}}%was scale=0.4
    \caption{Four low-dimensional simplices.}
    \label{fig: simplices in different dimensions}
\end{figure}

A \emph{simplex} is a geometric object\footnote{It is important to note that our geometric notion of simplices is different from the topological notion of a simplex, which is only concerned with the structure of a simplex up to homeomorphisms or the abstract/combinatorial notion of simplices, for which simplices are simply a collection of subsets of elements with the collection being closed under the subset relation.} that generalizes points ($d=0$), line segments ($d=1$) and triangles ($d=2$) to any dimension $d$, see Figure \ref{fig: simplices in different dimensions}. The classic characterization ``{three non-collinear points in the plane determine a triangle}" translates to ``{$n$ affinely independent points in $\mathbb{R}^{n-1}$ determine a simplex}" in this generalized setting. More precisely, a set of $n$ points\footnote{These points may also lie in an $(n-1)$-dimensional subspace of a larger-dimensional latent space.} $s_i\in\mathbb{R}^{n-1}$ such that for any $j$ the vectors $\lbrace s_i-s_j\rbrace_{i\neq j}$ are linearly independent, determines a simplex $S$ as their convex hull. Such points $s_i$ are indexed by $i\in\mathcal{N}$ and are called the \emph{vertices} of $S$. We will also use the notation $S=[s_1~\dots~s_n]$ to denote the $(n-1)\times n$ vertex matrix with columns equal to the vertex vectors. 
\\
We will mainly be interested in \emph{equivalence classes} of simplices, where two simplices are equivalent (congruent) if their vertex matrices satisfy
\begin{equation}\label{eq: simplex equivalence relation}
S\sim S' \Leftrightarrow \exists O,x:S'=OS+xu^T
\end{equation}
for some orthogonal matrix $O$, i.e. with $O^TO=I$, and vector $x\in\mathbb{R}^{n-1}$. In other words, this describes equivalence with respect to rotations, reflections and translations of the simplex which are all angle and distance-preserving\footnote{Rotations, reflections and translations are rigid transformations, which are isometries of Euclidean space.} . We will denote an equivalence class of simplices by $\mathcal{S}=\lbrace S':S'\sim S\rbrace$ with $S$ a specific representative of the equivalence class, and refer to the equivalence class $\mathcal{S}$ as a \emph{Simplex} (upper-case) and to a particular representative $S\in\mathcal{S}$ as a simplex (lower-case).
\\
In practice, we can represent a simplex $S$ by its so-called Gram matrix $S^TS$, i.e. with entries equal to the inner-product between vertex vectors. This Gram matrix is independent of rotations and reflections of the underlying simplex, but it does depend on translations. Thus, in order to define a unique Gram matrix for a Simplex $\mathcal{S}$ we fix a translation $S\rightarrow S(I-uu^T/n)$. Since the vector $Su/n=c$ is the \emph{centroid} of simplex $S$, the translated simplex equals $S-cu^T$ and its centroid coincides with the origin of $\mathbb{R}^{n-1}$. We will refer to this simplex as the \emph{centered simplex}. Having fixed a canonical translation, we can now define the \emph{canonical Gram matrix} of a Simplex as
$$
M \triangleq (I-uu^T/n)^TS^TS(I-uu^T/n)\text{~for any~}S\in\mathcal{S},
$$
which is a unique Gram matrix given the Simplex $\mathcal{S}$; if we want to refer to a specific Simplex, we will also write $M(\mathcal{S})$ for the canonical Gram matrix. Similarly, we define the \emph{canonical pseudoinverse Gram matrix} of a Simplex as the (Moore-Penrose) pseudoinverse of its canonical Gram matrix, and write $M^\dagger$. These matrix representations of Simplices have the following property:
\begin{proposition}\label{proposition: Gram matrix of a Simplices is PSD}
The canonical (pseudoinverse) Gram matrix of a Simplex satisfies properties (i)$_\sigma$-(iii)$_\sigma$. Conversely, any such matrix is the canonical (pseudoinverse) Gram matrix of some Simplex.
\end{proposition}
\textbf{Proof:} Let $S'$ be a representative simplex of a Simplex $\mathcal{S}$ and center this simplex as $S=S'(I-uu^T/n)$. By construction, the canonical Gram matrix $M=S^TS$ is now a symmetric, positive semidefinite matrix. Furthermore, the vertices of a simplex are affinely independent, which means that for any $j$ the set $\lbrace s_i-s_j\rbrace_{i\neq j}$ is linearly independent and thus that the matrix $S(I-e_ju^T)$ and also $(I-e_ju^T)^TM(I-e_ju^T)$ must have rank $(n-1)$. The Gram matrix will thus have $\operatorname{rank}(M)\geq (n-1)$. Finally, the product $Mu=(I-uu^T/n)S^{'T}S'(I-uu^T/n)u=0$ shows that $M$ has one zero eigenvalue with corresponding constant eigenvector. As a result, its pseudoinverse $M^\dagger$ satisfies the same properties.
\\
For the converse, let $M$ be any positive semidefinite matrix with a single zero eigenvalue and corresponding constant eigenvector, which can be decomposed as $M=\sum_{k=1}^{n-1}\mu_kz_kz_k^T$. From this decomposition, we find the Gram form $M=S^TS$ where $(S_i)_k=(z_k)_i\sqrt{\mu_k}$. The rows of $S$ are thus the (scaled) non-constant eigenvectors of $M$ such that $Sx=0$ if and only if $x$ is a constant (possibly zero) vector. Let $s_1,s_2,\dots,s_n$ be the columns of the matrix $S$; we show that for any $j\in\mathcal{N}$ the set of vectors $\lbrace s_i-s_j\rbrace_{i\neq j}$ is linearly independent. Assume that $\sum_{i\neq j} x_i(s_i-s_j)=0$. Then by letting $x_j=-\sum_{i\neq j}x_i$, we have $\sum_{i=1}^n x_is_i=0$ or equivalently $Sx=0$. Then $x$ must be a multiple of $u$. However, by construction $x$ is also orthogonal to $u$ and thus $x=0$ must hold, which proves the linear independence of $\lbrace s_i-s_j\rbrace_{i\neq j}$. The points $\lbrace s_i\rbrace_{i\in\mathcal{N}}$ are thus vertices of a simplex. Moreover, since $Su=0$ this simplex is centered and $M$ is thus the canonical Gram matrix of the Simplex $\mathcal{S}$ with representative $S$.
\\
Finally, if a matrix $M$ satisfies properties (i)$_\sigma$-(iii)$_\sigma$ then so will its pseudoinverse $M^\dagger$. By the previous derivation, we then know that $M^\dagger$ is the canonical Gram matrix of a Simplex $\mathcal{S}$ and thus that $M=(M^\dagger)^\dagger$ is the canonical pseudoinverse Gram matrix of $\mathcal{S}$.\hfill$\square$
\\
From Propositions \ref{proposition: equivalent characterizations for Laplacian matrices} and \ref{proposition: Gram matrix of a Simplices is PSD} it follows that every (pseudoinverse) Laplacian can be seen as the canonical Gram matrix of a Simplex. Conversely however, the canonical Gram matrix of a Simplex is only `somewhat like' a Laplacian matrix; more precisely it looks like a Laplacian with respect to the spectral properties (i)$_\sigma$-(iii)$_\sigma$ but can miss the sign property (ii) in general. The canonical pseudoinverse Gram matrix of a Simplex is another candidate Laplacian since it satisfies the spectral properties as well. In fact, we will show that the sign property of this pseudoinverse Gram matrix is related to angles in the Simplex. 
\\
A \emph{face} of a simplex $S$ is defined as the (sub)simplex determined by a subset $\mathcal{V}\subseteq\mathcal{N}$ of $v=\vert\mathcal{V}\vert$ vertices, and is denoted by $S_\mathcal{V}$. In Section \ref{section: faces of a simplex} we show that a face of a simplex is indeed again a simplex. Faces with $v=(n-1)$ vertices are called \emph{facets} and since these facets lie in a $(n-2)$-dimensional hyperplane of $\mathbb{R}^{n-1}$, they determine pairwise angles. In particular, we define the \emph{dihedral angle} $\phi_{ij}$ between facets $S_{\lbrace i\rbrace^c}$ and $S_{\lbrace j\rbrace^c}$ as the interior angle (with respect to the simplex) between these facets. Since all representative simplices of a Simplex are congruent, they have the same set of dihedral angles and we can define the dihedral angles of a Simplex as those of any representative simplex. We find the following relation between the canonical pseudoinverse matrix and dihedral angles:
\begin{property}\label{property: dihedral angles and pseudoinverse Gram}
The canonical pseudoinverse Gram matrix relates to the dihedral angles of a Simplex as
\begin{align*}
    (M^\dagger)_{ij}>0&\Leftrightarrow \phi_{ij}>\pi/2\textup{~is obtuse}.\\
    (M^\dagger)_{ij}=0&\Leftrightarrow\phi_{ij}=\pi/2\textup{~is right}.\\
    (M^\dagger)_{ij}<0&\Leftrightarrow\phi_{ij}<\pi/2\textup{~is acute}.
\end{align*}
\end{property}
\textbf{Proof:} Let $S$ be a centered representative simplex of Simplex $\mathcal{S}$, and $S^\dagger$ the Moore-Penrose pseudoinverse of its vertex matrix. Thus $S$ is an $(n-1)\times n$ vertex matrix with columns $S=[s_1~\dots~s_n]$ and $S^\dagger$ an $n\times (n-1)$ matrix and we write $S^\dagger=[\tilde{s}_1~\dots~\tilde{s}_n]^T$, i.e. $S^\dagger$ has rows $\tilde{s}_i^T$. By their pseudoinverse relation, these matrices satisfy $S^\dagger S=I-uu^T/n$. A facet $S_{\lbrace i\rbrace^c}$ of $S$ determines a hyperplane
$$
\mathcal{H}_i = \left\lbrace Sx:\forall x\in\mathbb{R}^n\text{~s.t.~}x_i=0\text{~and~}u^Tx=0\right\rbrace
$$
of vectors parallel with $S_{\lbrace i\rbrace^c}$. Indeed, any vector $x$ with $x_i=0$ and $x^Tu=0$ can be decomposed as $x=\alpha(y - y')$ for some scalar $\alpha\in\mathbb{R}$ and with vectors $y,y'$ with $y_i=y'_i=0$ and $u^Ty=u^Ty'=1$ such that $Sy,Sy'$ are points in the facet $S_{\lbrace i\rbrace^c}$. The vector $Sx=\alpha(Sy-Sy')$ is parallel to the line through these points and thus parallel to the facet.
\\
The rows of the pseudoinverse vertex matrix $S^\dagger$ satisfy 
$$
\tilde{s}_i^T(Sx) = e_i^T(S^{\dagger}S)x = e_i^Tx - e_i^Tuu^Tx/n = 0
$$
when $x_i=0$ and $u^Tx=0$, in other words $\tilde{s}_i\perp\mathcal{H}_i$. Since $\tilde{s}_i^Ts_i=1-1/n>0$ and the simplex $S$ is centered, this shows that $\tilde{s}_i$ is the \emph{inner-normal} vector of the facet $S_{\lbrace i\rbrace^c}$. As illustrated in the figure below, the dihedral (inner) angle $\phi_{ij}$ between a pair of facets can be calculated from the angle between these corresponding inner normal vectors as
$$
\cos(\pi-\phi_{ij}) = \frac{\tilde{s}_i^T\tilde{s}_j}{\Vert \tilde{s}_i\Vert\Vert \tilde{s}_j\Vert} = \frac{(M^\dagger)_{ij}}{\sqrt{(M^\dagger)_{ii}(M^{\dagger})_{jj}}},
$$
where the matrix product pseudoinverse satisfies $(S^\dagger S^{\dagger T}) = (S^TS)^\dagger$ because it is a product between transposed matrices $S$ and $S^T$, see \cite{Greville_pinv}. This
shows that the sign of $(M^\dagger)_{ij}$ determines the acute/right/obtuseness of dihedral angle $\phi_{ij}$.\hfill$\square$
\begin{figure}
    \centering\fbox{
    \includegraphics[scale=1.2]{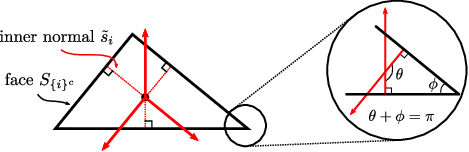}}
    \label{fig: inner normal vectors}
\end{figure}
\\
We stress that for Property \ref{property: dihedral angles and pseudoinverse Gram} to hold it is crucial that the canonical Gram matrix is based on a centered representative simplex, and thus that the canonical (pseudoinverse) Gram matrix has a zero eigenvalue corresponding to the constant eigenvector. In other words, the convenience of Property \ref{property: dihedral angles and pseudoinverse Gram} supports this specific choice of canonical Gram matrix.
%%%
%
%%%
\subsection{Simplices and Laplacians}\label{Section: simplices and laplacians}
Following Property \ref{property: dihedral angles and pseudoinverse Gram}, we know that if a Simplex is \emph{hyperacute}, i.e. when all of its dihedral angles are non-obtuse, then its canonical pseudoinverse Gram matrix will have non-positive off-diagonals. In other words, as the spectral properties (i)$_\sigma$-(iii)$_\sigma$ are automatically satisfied for a canonical pseudoinverse Gram matrix of a Simplex and thus a hyperacute Simplex in particular, we find that all properties of a Laplacian matrix are satisfied. We thus have:
\begin{lemma}\label{Lemma: Laplacians are Gram matrices}
The canonical pseudoinverse Gram matrix of every hyperacute Simplex is a Laplacian matrix. Conversely, any Laplacian matrix is the canonical pseudoinverse Gram matrix of a hyperacute Simplex.
\end{lemma}
\textbf{Proof:} By Proposition \ref{proposition: Gram matrix of a Simplices is PSD}, any canonical pseudoinverse Gram matrix satisfies properties (i)$_\sigma$-(iii)$_\sigma$. Moreover, for a hyperacute Simplex, the canonical pseudoinverse Gram matrix satisfies the sign property (ii) as well, such that this matrix is a Laplacian by Proposition \ref{proposition: equivalent characterizations for Laplacian matrices}. \\
Conversely, since the Laplacian matrix satisfies properties (i)$_\sigma$-(iii)$_\sigma$ it is a canonical pseudoinverse Gram matrix of a Simplex by Proposition \ref{proposition: Gram matrix of a Simplices is PSD}. Furthermore, the non-positive off-diagonal entries of the Laplacian imply by Property \ref{property: dihedral angles and pseudoinverse Gram} that this Simplex is hyperacute.\hfill$\square$
\\
Phrased differently, Lemma \ref{Lemma: Laplacians are Gram matrices} describes a correspondence between Laplacian matrices and hyperacute Simplices, which is best summarized as follows:
\begin{theorem}[Fiedler \cite{Fiedler_identity_and_equivalence}]\label{Theorem: bijection between Laplacians and Simplices}
There is a bijection between
\begin{enumerate}
    \item Laplacian matrices of $n$ dimensions
    \item hyperacute Simplices on $n$ vertices
\end{enumerate}
\end{theorem}
\textbf{Proof:} This follows from Lemma \ref{Lemma: Laplacians are Gram matrices}.\hfill$\square$
\\
The bijection described by Theorem \ref{Theorem: bijection between Laplacians and Simplices} is constructive in a straightforward way: for a given Simplex $\mathcal{S}$, we can always construct the corresponding Laplacian as the canonical pseudoinverse Gram matrix $Q=\big((I-uu^T/n)S^TS(I-uu^T/n)\big)^\dagger$ for any representative simplex $S$. For a given Laplacian $Q$ we can construct the corresponding Simplex as the equivalence class of the simplex with vertices $(s_i)_k=(z_k)_i\sqrt{1/\mu_k}$. Theorem \ref{Theorem: bijection between Laplacians and Simplices} allows us to speak, unambiguously, about the Laplacian of a Simplex and the Simplex of a Laplacian. 
\\
Theorem \ref{Theorem: bijection between Laplacians and Simplices} was discovered by Miroslav Fiedler in \cite{Fiedler_aggregation} (see also \cite{Fiedler_identity_and_equivalence}, \cite[Sec. 3.3]{Fiedler_book}) and sets up a rich connection between graph theory, linear algebra and (simplex) geometry, with many interesting implications described in \cite{Fiedler_book,krl_simplex}. In the next section we introduce the effective resistance, which is a key concept with valuable interpretations in graphs, Laplacians and simplices, as well as providing another perspective from which to understand Theorem \ref{Theorem: bijection between Laplacians and Simplices}.
%%%
%
%%%
\subsection{Effective resistances and Fiedler's identity}\label{section: effective resistance and fiedler}
The effective resistance was originally defined for resistive electrical circuits as the voltage measured between a pair of terminals $i$ and $j$ in the circuit, when a unit current is forced between these terminals. In other words, it captures the resistive effect of the whole network with respect to these terminals into a single `effective' resistance value (resistance = voltage/current) \cite{Biggs_potential}. For planar graphs, the effective resistance can also be obtained by applying a sequence of basic graph modifications which leave the effective resistance unchanged, until a single link remains between the terminals of choice (with resistance between these terminals equal to the effective resistance) \cite{Colin_de_verdiere, curtis, Chang_thesis}. From the perspective of random walks on the graph, the effective resistance can be calculated (up to a constant factor) as the average time it takes a random walker to go from one node to another, and back, the so-called commute time between these nodes \cite{Chandra}.
\\
Due to Kirchhoff's translation of the circuit equations in terms of the graph Laplacian, the effective resistance $\omega_{ij}$ between a pair of nodes $i$ and $j$ in a graph with Laplacian $Q$ can be calculated, and thus defined, as follows \cite{Ellens, Bapat_inverse_of_weighted_graphs, Ghosh}:
\begin{equation}\label{eq: definition effective resistance}
\omega_{ij} = (e_i-e_j)^TQ^\dagger(e_i-e_j),
\end{equation}
with basis vectors $(e_i)_k=1$ if $k=i$ and zero otherwise. In other words, the effective resistance can be found from a quadratic product with the pseudoinverse Laplacian. The $n\times n$ matrix $\Omega$ containing all effective resistances as its entries $(\Omega)_{ij}=\omega_{ij}$, is called the \emph{resistance matrix}.
\\
An important property of the effective resistance is that it provides another `bridge' between Laplacian matrices and Simplices. If we introduce the Gram representation $Q^\dagger = S^TS$ of the pseudoinverse Laplacian in definition \eqref{eq: definition effective resistance}, we find that
\begin{equation}\label{eq: effective resistance as edge lengths}
\omega_{ij} =(e_i-e_j)^TQ^\dagger(e_i-e_j) = \Vert S(e_i-e_j)\Vert^2=\Vert s_i-s_j\Vert^2.
\end{equation}
In other words, for a graph $G$ with Laplacian matrix $Q$ and corresponding Simplex $\mathcal{S}$, the effective resistance between a pair of nodes in $G$ equals the squared distance between the corresponding pair of vertices in $\mathcal{S}$; equivalently, the vertices of $\mathcal{S}$ are an embedding of the nodes of $G$ into $\mathbb{R}^{n-1}$ where the effective resistance matrix thus plays the role of the \emph{squared Euclidean distance matrix} of the simplex $\mathcal{S}$. In fact, since distances between vertices are invariant with respect to reflections, rotations and translations, and contain all information necessary to reconstruct a Simplex, the resistance matrix $\Omega$ characterizes the equivalence class as
$$
\mathcal{S} = \left\lbrace S' : \Vert S'(e_i-e_j)\Vert^2=\omega_{ij} \text{~for all~}i,j \right\rbrace.
$$
We will also write $\mathcal{S}(\Omega)$ if we want to further specify the effective resistance matrix and conversely for $\Omega(\mathcal{S})$. 
\\
The effective resistance allows the bijection between simplices, graphs and Laplacian matrices to be summarized beautifully by the following identity:
\begin{theorem}[Fiedler's identity]\label{Theorem: Fiedlers identity}
For a weighted graph $G$ with Laplacian matrix $Q$ and Simplex $\mathcal{S}$ with resistance matrix $\Omega$, the following identity holds
\begin{equation}\label{eq: fiedlers identity}
-\frac{1}{2}\begin{pmatrix}0&u^T\\u&\Omega\end{pmatrix} = \begin{pmatrix}4R^2&-2r^T\\-2r&Q\end{pmatrix}^{-1}
\end{equation}
where $r=\tfrac{1}{2}Q\zeta+\tfrac{u}{n}$ with $\zeta=\operatorname{diag}(Q^\dagger)$ determines the \emph{circumcenter} of $\mathcal{S}$ as $Sr$, and with $R=\sqrt{\tfrac{1}{4}\zeta^TQ\zeta + \tfrac{u^T\zeta}{n}}$ the \emph{circumradius} of $\mathcal{S}$.\\
\end{theorem}
\textbf{Proof:} We will conduct the proof in four steps.\\ 
\textbf{Step 1} We start by showing that $\Omega$ is invertible: from definition \eqref{eq: definition effective resistance} of effective resistances, the resistance matrix can be decomposed as 
\begin{equation}\label{eq: resistance matrix decomposition}
\Omega = u\zeta^T + \zeta u^T -2Q^\dagger\text{~with $\zeta=\operatorname{diag}(Q^\dagger)$}.
\end{equation}
From this decomposition we find that $x^T\Omega x = -2x^TQ^\dagger x< 0$ for all vectors $x^Tu=0$, where the inequality follows from the spectral properties of the pseudoinverse Laplacian. Furthermore, from the fact that the resistance matrix has positive off-diagonal entries we know that $u^T\Omega u>0$. Combining these inequalities shows that the resistance matrix has precisely one positive eigenvalue and $n-1$ negative eigenvalues (i.e. it is an elliptic matrix) and is thus non-singular.
\\
\textbf{Step 2} Following the invertibility of the resistance matrix, we can solve the equation $\Omega x= u$ for $x$: using decomposition \eqref{eq: resistance matrix decomposition} of the resistance matrix, we find that
\begin{align*}
&\hphantom{\Leftrightarrow} u = (u\zeta^T+\zeta u^T-2Q^\dagger)x  
\\
&\Leftrightarrow u = u(\zeta^Tx)+\left(I-\frac{uu^T}{n}+\frac{uu^T}{n}\right)\zeta (u^Tx)-2Q^\dagger x
\\
&\Leftrightarrow u\left(1 - \zeta^Tx - \frac{(u^T\zeta)(u^Tx)}{n}\right) = \left(I-\frac{uu^T}{n}\right)\zeta (u^Tx) - 2Q^\dagger x
\end{align*}
Since the lefthandside of this equation is a multiple of $u$ while the righthandside is perpendicular to $u$, it follows that both sides must be zero, and thus
\begin{equation}
\begin{dcases}
\zeta^Tx = 1-\frac{(u^T\zeta)(u^Tx)}{n} \label{eq: equation for zeta^Tx}
\\
2Q^\dagger x = \left(I-\frac{uu^T}{n}\right)\zeta (u^Tx)
\end{dcases}
\end{equation}
Multiplying both sides of the second expression in \eqref{eq: equation for zeta^Tx} by $Q$, we find
\begin{equation}\label{eq: expression for x}
2\left(I-\frac{uu^T}{n}\right)x = Q\zeta(u^Tx)\Leftrightarrow x = (u^Tx)(\tfrac{1}{2}Q\zeta+\tfrac{u}{n}).
\end{equation}
After multiplication by $\zeta^T$ (which is not the zero vector, see next step) and introducing the value for $\zeta^Tx$ from the first expression in \eqref{eq: equation for zeta^Tx}, this yields
\begin{equation}\label{eq: expression for u^Tx}
(u^Tx)\left(\frac{1}{2}\zeta^TQ\zeta+2\frac{u^T\zeta}{n}\right) = 1.
\end{equation}
Introducing $r=\tfrac{1}{2}Q\zeta + \tfrac{u}{n}$ into expression \eqref{eq: expression for x} for $x$, and $2R^2=\frac{1}{2}\zeta^TQ\zeta+2\tfrac{u^T\zeta}{n}$ into expression \eqref{eq: expression for u^Tx} for $(u^Tx)$, this becomes
$$
x = (u^Tx)r \text{~and~}(u^Tx)2R^2 = 1.
$$
If we assume that $0<2R^2$ (which is confirmed at the end of this step), these expressions imply that $x=\frac{r}{2R^2}$ and thus give the equation
\begin{equation}\label{eq: relation between r and R}
\Omega r = 2R^2 u,
\end{equation}
and the relation $2R^2= r^T\Omega r$. By inverting $\Omega$, multiplying by $u^T$ and invoking the unit-sum property $u^Tr=1$ (as follows from $Qu=0$), equation \eqref{eq: relation between r and R} yields an alternative definition of $R,r$ in terms of the resistance matrix:
$$
2R^2 = \frac{1}{u^T\Omega^{-1}u}\quad\text{~and~}\quad r = \frac{\Omega^{-1}u}{u^T\Omega^{-1}u}.
$$
\\
We now verify the bound for $2R^2=\frac{1}{2}\zeta^TQ\zeta+2\tfrac{u^T\zeta}{n}$: by positive semidefiniteness of $Q$ and rewriting $u^T\zeta = \tr(Q^\dagger)$ -- which is a strictly positive trace by the spectral properties of $Q^\dagger$ and $n\geq 2$ -- we find that $2R^2>0$ as required.
\\
\textbf{Step 3} Next, we show that the matrix $\left(\begin{smallmatrix}0&u^T\\u&\Omega\end{smallmatrix}\right)$ is invertible: assuming otherwise, there must exist a scalar $\alpha$ and vector $y$, not both zero, such that
\begin{equation}\label{eq: omega extended is invertible}
\begin{pmatrix}
0&u^T\\u&\Omega
\end{pmatrix}\begin{pmatrix}
-\alpha\\ y
\end{pmatrix} = \begin{pmatrix}
0\\0
\end{pmatrix}\Leftrightarrow\begin{cases}
u^Ty=0\\
\Omega y = \alpha u,
\end{cases}
\end{equation}
First, assuming $y\neq 0$ then by non-singularity of $\Omega$ and the second equation $\Omega y = \alpha u$ in \eqref{eq: omega extended is invertible}, we know that $\alpha\neq 0$ must hold as well. But from equation \eqref{eq: relation between r and R} we then find that $y=\alpha\Omega^{-1}u=\alpha r/(2R^2)$ which means that $u^Ty=\alpha (u^Tr)/(2R^2)\neq 0$ since $u^Tr=1, \alpha\neq 0$ and $0<2R^2<\infty$. This is in contradiction with $u^Ty=0$ in equation \eqref{eq: omega extended is invertible}, hence $y\neq 0$ is not possible. But if $y=0$ then by equation \eqref{eq: omega extended is invertible} also $\alpha=0$ must hold, in contradiction with the assumption that $y$ and $\alpha$ are not both zero. It thus follows that $\left(\begin{smallmatrix}0&u^T\\u&\Omega\end{smallmatrix}\right)$ is invertible.
\\
\textbf{Step 4} Finally, we can verify the proposed matrix inverse \eqref{eq: fiedlers identity}: combining expression \eqref{eq: relation between r and R} as $-2R^2u+\Omega r=0$ and the unit-sum property $u^Tr=1$ into a single matrix expression, we find
\begin{align*}
\begin{cases}
(-2R^2).0 + u^Tr = 1\\
-2R^2u + \Omega r = 0
\end{cases}&\Leftrightarrow 
\begin{pmatrix}
0&u^T\\u&\Omega
\end{pmatrix}\begin{pmatrix}
-2R^2\\ r
\end{pmatrix}=\begin{pmatrix}1\\0\end{pmatrix}
\\
&\Leftrightarrow \begin{pmatrix}0& u^T\\ u & \Omega\end{pmatrix}^{-1}\begin{pmatrix}1\\0\end{pmatrix} = \begin{pmatrix}-2R^2\\r\end{pmatrix}
\\
&\Leftrightarrow \begin{pmatrix}0 & u^T\\ u & \Omega\end{pmatrix}^{-1} = \begin{pmatrix}-2R^2 & *\\r & *\end{pmatrix}\\
&\Leftrightarrow \begin{pmatrix}0 & u^T\\ u & \Omega\end{pmatrix}^{-1} = \begin{pmatrix}-2R^2 & r^T\\r & A\end{pmatrix}
\end{align*}
where the last step follows by symmetry, and where $A$ is some symmetric matrix that is yet to be determined. From the matrix product $\left(\begin{smallmatrix}0&u^T\\u&\Omega\end{smallmatrix}\right)\left(\begin{smallmatrix}-2R^2&r^T\\r& A \end{smallmatrix}\right)=I$ we then find that $A$ must satisfy
\begin{equation}\label{eq: relation between Omega and M}
ur^T+\Omega A =  I\text{~and~} Au = 0.
\end{equation}
Left-multiplying the first equation by $Q$ and making use of $Qu=Au=0$ and \eqref{eq: resistance matrix decomposition} retrieves the Laplacian matrix $A=-\tfrac{1}{2}Q$. Since the inverse is uniquely determined this verifies Fiedler's identity \eqref{eq: fiedlers identity}. To conclude, the interpretation of $R$ and $Sr$ as the respective circumradius and circumcenter of $\mathcal{S}$ are proven by Fiedler \cite[Cor. 1.4.13]{Fiedler_book}, and in \cite{krl_covariance}.
\\
We remark that while deriving Theorem \ref{Theorem: Fiedlers identity} as above requires some work, it is much easier to verify once the expressions for $r$ and $R$ are known; multiplying both sides of Fiedler's identity \eqref{eq: fiedlers identity} retrieves the identity matrix and thus constitutes a more direct but less transparent proof.\hfill$\square$
\\
While $G,Q,\mathcal{S},\Omega$ are mutually interchangeable, it is somehow natural to think of $G$ and $\mathcal{S}$ as the combinatorial and geometric objects of interest, with $Q$ and $\Omega$ as convenient and practical algebraic representations. From this perspective, Fiedler's identity sets up the \emph{graph-Simplex correspondence} via a direct inverse relation between their respective representations. Importantly, this direct algebraic identity also gives a way to characterize how certain operations on Laplacian matrices translate to operations on the resistance matrix and vice versa. This fact will be key in our discussion of the Schur complement in Section \ref{Section: faces and Schur} and plays a major role in the final structure underlying the proof in Section \ref{Section: geometric proof of resistance distance}.
\\
Apart from providing a concise summary of the equivalences, Fiedler's identity brings up a number of additional interesting results. It show how the circumradius and circumcenter of $\mathcal{S}$  are expressed in terms of the (pseudoinverse) Laplacian $Q$, and introduces the matrix $\left(\begin{smallmatrix}0&u^T\\u&\Omega\\\end{smallmatrix}\right)$ which is known as the \emph{Cayley-Menger} matrix and is related to the volume of $\mathcal{S}$ \cite{Fiedler_book}. We furthermore believe that the vector $r$ and scalar $R$ are important algebraic objects associated to a graph, worthy of a deeper study, e.g. as initiated in \cite{krl_covariance}.
\\
\textit{Remark:} While stated in terms of Laplacian and resistance matrices, Fiedler's identity also holds for non-hyperacute Simplices since the proof of Theorem \ref{Theorem: Fiedlers identity} does not rely on the sign of Laplacian entries or the dihedral angles of a Simplex. For a non-hyperacute Simplex $\mathcal{S}$, this yields the matrix identity 
\begin{equation}\label{eq: fiedlers identity general simplices}
-\frac{1}{2}\begin{pmatrix}0&u^T\\u&D(\mathcal{S})\end{pmatrix} = \begin{pmatrix}4R^2&-2r^T\\-2r&M^\dagger(\mathcal{S})\end{pmatrix}^{-1}
\end{equation}
between the squared Euclidean distance matrix $D(\mathcal{S})$ and canonical pseudoinverse Gram matrix $M^\dagger(\mathcal{S})$ of $\mathcal{S}$. More generally, for any invertible matrix $A$ with $u^T A^{-1}u\neq 0$ we find that the matrix $\left(\begin{smallmatrix}0 &u^T\\u&A \end{smallmatrix}\right)$ is invertible; consequently, certain results and techniques that follow from Fiedler's identity will also be applicable to the analysis of $A$. In particular, we believe that this might be relevant to the theory of \emph{magnitude of metric spaces} \cite{Leinster}, where for matrices of the form $(A_t)_{ij} = \operatorname{exp}(-d(i,j)t)$ for some $t>0$ and metric $d$, the magnitude is defined as $\vert A_t\vert :=u^TA_t^{-1}u$ if this inverse exists.
\\
\textit{Remark:} The matrix inverse of resistance matrices was rediscovered later, independent of Fiedler's work, by Graham and Lov{\'{a}}sz \cite{Graham_Lovasz_inverse_resistance_tree} for tree graphs and by Bapat \cite{Bapat_inverse_of_weighted_graphs} for general weighted graphs. They describe the elegant relation $\Omega^{-1}=-\frac{1}{2}Q + \frac{rr^T}{2R^2}$ which follows from expression \eqref{eq: relation between Omega and M} in the proof of Theorem \ref{Theorem: Fiedlers identity} by left-multiplication with $\Omega^{-1}$ and from \eqref{eq: relation between r and R}, or more directly from Fiedler's identity by taking the Schur complement introduced in Section \ref{Section: Schur complement}.
%%%
%
%%%
\section{Maps between Laplacians, maps between Simplices}\label{Section: faces and Schur}
In the previous section we defined graphs, Laplacians, resistance matrices and simplices and the relations between these different objects. Here, we will study instead the relation between instances of the same type of object, e.g. between pairs of Laplacian matrices and pairs of Simplices. Starting from the `face relation' between simplices, we find a corresponding `submatrix' relation for Laplacian matrices (the Schur complement) which retains important properties of the original Laplacian.
%%%
%
%%%
\subsection{Faces of a Simplex}\label{section: faces of a simplex}
We recall that the face of a simplex $S$ is the convex hull of a subset $\mathcal{V}$ of the vertices of $S$, and is denoted by $S_\mathcal{V}$. Similarly, we define the \emph{face of a Simplex} $\mathcal{S}/\mathcal{V}^c$ as the equivalence class
\begin{equation}\label{eq: face equivalence definition}
\mathcal{S}/\mathcal{V}^c \triangleq \lbrace S':S'\sim S_\mathcal{V}\rbrace\text{~for any~}S\in\mathcal{S}
\end{equation}
with the equivalence $\sim$ as defined in \eqref{eq: simplex equivalence relation}. We will also call $\mathcal{S}/\mathcal{V}^c$ the $\mathcal{V}$-face of $\mathcal{S}$. The choice for this notation will be explained later. An important property of faces is the following:
\begin{property}\label{Property: face of a Simplex is Simplex}
The face of a Simplex is again a Simplex.
\end{property}
\textbf{Proof:} Let $S=[s_1~\dots~s_n]$ be a representative simplex of $\mathcal{S}$. Since $S$ is a simplex, the vertices are affinely independent, i.e. $\lbrace s_i-s_j\rbrace_{i\neq j}$ is a linearly independent set for all $j$, and $\operatorname{rank}(S)=n-1$. Removing any subset $\mathcal{V}^c$ of vectors from this set of vector differences yields the set $\lbrace s_i-s_j\rbrace_{i\in\mathcal{V}\backslash j}$ whose rank is reduced by at most $(n-v)$ (by properties of the rank), with $v=\vert\mathcal{V}\vert$. As a consequence, this set has rank $\geq (v-1)$ where equality must hold since the rank cannot exceed the cardinality of the set. Since this holds for any $j$, $S_\mathcal{V}$ is again a simplex and the equivalence class $\lbrace S':S'\sim S_\mathcal{V}\rbrace$ is a Simplex. \hfill$\square$ 
\\
This property is a fundamental characteristic of simplices; in fact, `closure under taking subsets' is part of the axiomatic definition of a simplex in the study of abstract simplices and simplicial complexes. An example of faces of a Simplex is shown in Figure \ref{fig: graph simplex and subsets}. \\
Definition \eqref{eq: face equivalence definition} for faces follows from the perspective of simplices as being specified by a given set of vertices. An alternative definition of the face of a Simplex follows from specifying the squared Euclidean distance matrix $D$ between the vertices instead. A face determined by a subset $\mathcal{V}$ of the vertices then corresponds to a submatrix of the distance matrix
\begin{equation}\label{eq: face distance definition}
D(\mathcal{S}/\mathcal{V}^c)=(D(\mathcal{S}))_{\mathcal{V}\mathcal{V}},
\end{equation}
and we have $\mathcal{S}/\mathcal{V}^c = \mathcal{S}((D)_{\mathcal{V}\mathcal{V}})$. This alternative description based on taking submatrices of the distance matrix clearly highlights the following property of the face relation:
\begin{property}[composition]\label{Property: composition property of faces}
The face relation between Simplices can be composed: for any $\mathcal{W}\subseteq\mathcal{V}\subseteq\mathcal{N}$, the $\mathcal{W}$-face of the $\mathcal{V}$-face of a Simplex $\mathcal{S}$ equals the $\mathcal{W}$-face of $\mathcal{S}$; in other words 
$$
\mathcal{S}/(\mathcal{N}\backslash\mathcal{W}) = [\mathcal{S}/(\mathcal{N}\backslash\mathcal{V})]/(\mathcal{V}\backslash\mathcal{W}).
$$
\end{property}
\textbf{Proof:} Let $\lbrace s_i\rbrace_{i\in\mathcal{N}}$ be the vertices of a representative simplex of $\mathcal{S}$. Following definition \eqref{eq: face equivalence definition} of faces of a Simplex, the $\mathcal{W}$-face of the $\mathcal{V}$-face of $\mathcal{S}$ is represented by a simplex with vertices $\lbrace s_i\rbrace_{i\in\mathcal{W}\subseteq\mathcal{V}}$. This is just the simplex with vertices $\lbrace s_i\rbrace_{i\in\mathcal{W}}$ which by definition \eqref{eq: face equivalence definition} is a representative simplex of the $\mathcal{W}$-face of $\mathcal{S}$. Since every representative simplex corresponds to a unique Simplex (equivalence classes determine a partition), both Simplex faces are equal as required.\hfill$\square$
\\
A property of this type is also called a \emph{quotient property} which supports our notation of faces as a `quotient' over a subset of the vertices.
\\
We remark that in terms of the distance matrix of a simplex, the composition property of faces is reflected in the submatrix relation
\begin{equation}\label{eq: distance submatrix for Simplex face}
(D(\mathcal{S}))_{\mathcal{W}\mathcal{W}} = [(D(\mathcal{S}))_{\mathcal{V}\mathcal{V}}]_{\mathcal{W}\mathcal{W}}.
\end{equation}
%%%
%
%%%
\subsection{The Schur complement}\label{Section: Schur complement}
In this section, we show how the two definitions of faces, via a subset of the vertices or via a submatrix of the distance matrix, lead to complementary expressions for the so-called Schur complement of a (Laplacian) matrix. From definition \eqref{eq: face equivalence definition}, a first expression for the canonical Gram matrix of a face follows:
\begin{proposition}\label{proposition: canonical Gram matrix of a face}
The canonical Gram matrix of a face of a Simplex $\mathcal{S}$ is equal to
\begin{equation}\label{eq: Gram matrix of face}
M(\mathcal{S}/\mathcal{V}^c) = (I-uu^T/v)(M(\mathcal{S}))_{\mathcal{V}\mathcal{V}}(I-uu^T/v).
\end{equation}
\end{proposition}
\textbf{Proof:} Let $S$ be the vertex matrix of a representative of the Simplex $\mathcal{S}$. Restricting to the columns corresponding to the vertices in the face, we get the vertex matrix $(S)_{*\mathcal{V}}$ of $\mathcal{S}/\mathcal{V}^c$ which has all rows of $S$ (denoted by the subscript $*$) and only those columns in $\mathcal{V}$. The canonical Gram matrix \eqref{eq: Gram matrix of face} of the face is then found by centering this vertex matrix as $(S)_{*\mathcal{V}}(I-uu^T/v)$ and using the fact that $(S)_{*\mathcal{V}}^T(S)_{*\mathcal{V}}=(M(\mathcal{S}))_{\mathcal{V}\mathcal{V}}$. \hfill$\square$
\\
As a consequence of \eqref{eq: Gram matrix of face}, we find that quadratic products with the canonical Gram matrix of a face correspond to quadratic products with the canonical Gram matrix of the Simplex as
\begin{equation}\label{eq: quadratic product Gram face}
x^TM(\mathcal{S}/\mathcal{V}^c)x = x^TM(\mathcal{S})x
\end{equation}
for all $x\in\mathbb{R}^n$ with $x^Tu=0$ and $x_i= 0$ if $i\notin\mathcal{V}$. This expression recovers the fact expressed by \eqref{eq: face distance definition} that the distances between vertices of a face are equal to the distances between the corresponding vertices of the Simplex -- which is obtained by choosing $x=(e_i-e_j)$ in \eqref{eq: quadratic product Gram face} to yield $x^TM(\mathcal{S})x = \Vert s_i-s_j\Vert^2$.
Proposition \ref{proposition: canonical Gram matrix of a face} furthermore says that the canonical pseudoinverse Gram matrix of a Simplex face equals 
$$
M^\dagger(\mathcal{S}/\mathcal{V}^c) = \left[\left(I-uu^T/v\right)(M(\mathcal{S}))_{\mathcal{V}\mathcal{V}}\left(I-uu^T/v\right)\right]^\dagger.
$$ In terms of algebraic operations, this corresponds to first taking a submatrix ($*_{\mathcal{V}\mathcal{V}}$) of the canonical Gram matrix followed by taking the pseudoinverse ($*^\dagger$). Performing these operations in reverse order will give rise to a second expression for the canonical pseudoinverse Gram matrix.
\\~\\
When combining matrix inverses and submatrices, the concept of Schur complements is relevant. For an invertible matrix $A=\left(\begin{smallmatrix}A_{\mathcal{V}\mathcal{V}} & A_{\mathcal{V}\mathcal{V}^c}\\A_{\mathcal{V}^c\mathcal{V}} & A_{\mathcal{V}^c\mathcal{V}^c}\end{smallmatrix}\right)$, the submatrix of its inverse equals \cite[Thm. 1.2]{Zhang_Schur_Complement}
\begin{equation}\label{eq: definition Schur complement invertible matrix}
(A^{-1})_{\mathcal{V}\mathcal{V}}  = \left(A_{\mathcal{V}\mathcal{V}} - A_{\mathcal{V}\mathcal{V}^c}(A_{\mathcal{V}^c\mathcal{V}^c})^{-1}A_{\mathcal{V}^c\mathcal{V}}\right)^{-1} \triangleq (A/\mathcal{V}^c)^{-1},
\end{equation}
where the introduced matrix\footnote{This Schur complement $A/\mathcal{V}^c$ is sometimes denoted by $A/A_{\mathcal{V}^c\mathcal{V}^c}$ instead.} $A/\mathcal{V}^c$ is called the \emph{Schur complement} of $A$ with respect to (the index subset) $\mathcal{V}$. The Schur complement is defined more generally for any matrix $A$ and subset $\mathcal{V}$ such that $A_{\mathcal{V}^c\mathcal{V}^c}$ is invertible. In the case of canonical pseudoinverse Gram matrices (and thus Laplacians), the spectral properties guarantee that $(M^\dagger)_{\mathcal{V}^c\mathcal{V}^c}$ is invertible\footnote{Assume for contradiction that $(M^\dagger)_{\mathcal{V}^c\mathcal{V}^c}$ is singular, then $(M^\dagger)_{\mathcal{V}^c\mathcal{V}^c}x=0$ holds for some vector $x$. But this implies that $M^\dagger \left(\begin{smallmatrix}0\\x\end{smallmatrix}\right)=0$ which contradicts the spectral properties (ii)$_\sigma$-(iii)$_\sigma$ of $M^\dagger$.} and thus that the Schur complement $M^\dagger/\mathcal{V}^c$, exist for any (nonempty) subset $\mathcal{V}\subset\mathcal{N}$. The Schur complement is a widely studied matrix operation in linear algebra and our discussion here will be limited to a number of properties which are relevant to our problem. For a general overview of the history and algebraic properties of the Schur complement, we refer to the excellent survey \cite{Zhang_Schur_Complement}.
\\
We now follow a second approach to identify the canonical pseudoinverse Gram matrix of a Simplex face. Combining the submatrix formula \eqref{eq: face distance definition} for the face distance matrix and Fiedler's identity for Simplices \eqref{eq: fiedlers identity general simplices}, we find
$$
-\frac{1}{2}\begin{pmatrix}
0 & u^T\\ u & D(\mathcal{S}/\mathcal{V}^c)
\end{pmatrix}
\overset{\eqref{eq: face distance definition}}{=}
-\frac{1}{2}\left[\begin{pmatrix}
0 & u^T\\ u & D(\mathcal{S})
\end{pmatrix}\right]_{\mathcal{V}_+\mathcal{V}_+}
\overset{\eqref{eq: fiedlers identity general simplices}}{=}
\left[\begin{pmatrix}
4R^2 & -2r^T\\
-2r & M^\dagger(\mathcal{S})
\end{pmatrix}^{-1}\right]_{\mathcal{V}_+\mathcal{V}_+}
$$
where $\mathcal{V}_+$ is the set of indices $\mathcal{V}$ together with the first row/column index. Using the Schur complement and invoking Fiedler's identity for the Simplex face then yields the following result:
\begin{proposition}\label{Proposition: Gram matrix of face is Schur complement}
The canonical pseudoinverse Gram matrix of a face of a Simplex is equal to the Schur complement of the canonical pseudoinverse Gram matrix of the Simplex; in other words:
\begin{equation}\label{eq: Gram matrix of face is Schur complement}
M^\dagger(\mathcal{S}/\mathcal{V}^c)=M^\dagger(\mathcal{S})/\mathcal{V}^c
\end{equation}
\end{proposition}
\textbf{Proof:} Starting from Fiedler's identity for a Simplex face we can derive
\begin{align*}
\begin{pmatrix}
4\tilde{R}^2 & -2\tilde{r}^T\\
-2\tilde{r} & M^\dagger(\mathcal{S}/\mathcal{V}^c)
\end{pmatrix}^{-1} 
&\overset{\eqref{eq: fiedlers identity general simplices}}{=}
-\frac{1}{2}\begin{pmatrix}
0 & u^T\\ u & D(\mathcal{S}/\mathcal{V}^c)
\end{pmatrix}
\\
&\overset{\eqref{eq: face distance definition}}{=}
\left[-\frac{1}{2}\begin{pmatrix}
0 & u^T\\ u & D(\mathcal{S})
\end{pmatrix}\right]_{\mathcal{V}_+\mathcal{V}_+}
\\
&\overset{\eqref{eq: fiedlers identity general simplices}}{=}
\left[\begin{pmatrix}
4R^2 & -2r^T\\
-2r & M^\dagger(\mathcal{S})
\end{pmatrix}^{-1}\right]_{\mathcal{V}_+\mathcal{V}_+}
\end{align*}
where $\tilde{R},\tilde{r}$ and $R,r$ are the circumradius and circumcenter coordinates of Simplex $\mathcal{S}/\mathcal{V}^c$ and $\mathcal{S}$, respectively\footnote{The specific values of $\tilde{R},R,\tilde{r},r$ are not important in this proof, see \cite{krl_covariance} for further details on how they are related.} (as in Theorem \ref{Theorem: Fiedlers identity}). Next, invoking the Schur complement definition \eqref{eq: definition Schur complement invertible matrix} for the $\mathcal{V}_+$ submatrix of the inverse, we find
\begin{align*}
&\begin{pmatrix}
4\tilde{R}^2 & -2\tilde{r}^T\\
-2\tilde{r} & M^\dagger(\mathcal{S}/\mathcal{V}^c)
\end{pmatrix}^{-1} 
\\
&=
\left[\begin{pmatrix}
4R^2 & -2r_{\mathcal{V}}^T\\-2r_{\mathcal{V}} & (M^\dagger(\mathcal{S}))_{\mathcal{V}\mathcal{V}}
\end{pmatrix} - 
\begin{pmatrix}-2r_{\mathcal{V}^c}^T\\ (M^\dagger(\mathcal{S}))_{\mathcal{V}\mathcal{V}^c}\end{pmatrix}[(M^\dagger(\mathcal{S}))_{\mathcal{V}^c\mathcal{V}^c}]^{-1}\begin{pmatrix}-2r_{\mathcal{V}^c} & (M^\dagger(\mathcal{S}))_{\mathcal{V}^c\mathcal{V}}\end{pmatrix}\right]^{-1}
\end{align*}
which by inverting both sides and considering the $\mathcal{V}$ submatrix yields
\begin{align*}
M^\dagger(\mathcal{S}/\mathcal{V}^c) &= (M^\dagger(\mathcal{S}))_{\mathcal{V}\mathcal{V}} - (M^\dagger(\mathcal{S}))_{\mathcal{V}\mathcal{V}^c}[(M^\dagger(\mathcal{S}))_{\mathcal{V}^c\mathcal{V}^c}]^{-1}(M^\dagger(\mathcal{S}))_{\mathcal{V}^c\mathcal{V}} 
\\
&= M^\dagger(\mathcal{S})/\mathcal{V}^c
\end{align*}
as required.\hfill$\square$
\\
As a result of Proposition \ref{Proposition: Gram matrix of face is Schur complement} and expression \eqref{eq: Gram matrix of face} we moreover have the following (known) alternative expression for the Schur complement of the canonical pseudoinverse Gram matrix:
\begin{equation}\label{eq: definition Schur complement as submatrix}
M^\dagger(\mathcal{S})/\mathcal{V}^c = M^\dagger(\mathcal{S}/\mathcal{V}^c) =  \left[(I-uu^T/v)(M(\mathcal{S}))_{\mathcal{V}\mathcal{V}}(I-uu^T/v)\right]^\dagger.
\end{equation}
This expression is complementary to definition \eqref{eq: definition Schur complement invertible matrix} as some important properties of the Schur complement are more apparent in the former expression than the latter.
\\
Next, from the relation between the Schur complement and the Simplex face relation we find the following composition property: 
\begin{property}[composition]\label{Property: composition property of Schur complement}
For any canonical pseudoinverse Gram matrix $M^\dagger$ and index subsets $\mathcal{W}\subseteq\mathcal{V}\subseteq\mathcal{N}$, the Schur complement of $M^\dagger/\mathcal{V}^c$ with respect to $\mathcal{W}$ is equal to the Schur complement of $M^\dagger$ with respect to $\mathcal{W}$; in other words, the Schur complement composes as
$$
M^\dagger/(\mathcal{N}\backslash\mathcal{W}) = [M^\dagger/(\mathcal{N}\backslash\mathcal{V})]/(\mathcal{V}\backslash\mathcal{W}).
$$
\end{property}
\textbf{Proof:} This composition property for the Schur complement of canonical pseudoinverse Gram matrices follows from the fact that these Schur complements correspond to faces of the Simplex (Proposition \ref{Proposition: Gram matrix of face is Schur complement}) together with the composition property of faces of a Simplex (Property \ref{Property: composition property of faces}). Repeated application of these two properties yields:
\begin{align*}
(M^\dagger(\mathcal{S}))/(\mathcal{N}\backslash\mathcal{W}) &\overset{Propos. \ref{Proposition: Gram matrix of face is Schur complement}}{=} M^\dagger(\mathcal{S}/(\mathcal{N}\backslash\mathcal{W}))\\
&\overset{Proper. \ref{Property: composition property of faces}}{=} M^\dagger\big([\mathcal{S}/(\mathcal{N}\backslash\mathcal{V})]/(\mathcal{V}\backslash\mathcal{W})\big)
\\
&\overset{Propos. \ref{Proposition: Gram matrix of face is Schur complement}}{=} \big[M^\dagger(\mathcal{S}/(\mathcal{N}\backslash\mathcal{V}))\big]/(\mathcal{V}\backslash\mathcal{W})
\\
&\overset{Propos. \ref{Proposition: Gram matrix of face is Schur complement}}{=} [(M^\dagger(\mathcal{S}))/(\mathcal{N}\backslash\mathcal{V})]/(\mathcal{V}\backslash\mathcal{W})
\end{align*}
as required.\hfill$\square$
\\
Property \ref{Property: composition property of Schur complement} holds in general for the Schur complement and was discovered by Emilie Haynesworth in \cite{Haynesworth}, where it was coined the \emph{quotient property} of the Schur complement and motivated the quotient notation of the Schur complement. One consequence of the composition property is that it allows the Schur complement $M^\dagger/\mathcal{V}^c$ with respect to any set $\mathcal{V}$ to be decomposed as a repeated application of the Schur complement with respect to complements of single indices as
\begin{equation}\label{eq: incremental Schur complement}
M^\dagger/\mathcal{V}^c = M^\dagger/\lbrace v_1\rbrace/\dots/\lbrace v_k\rbrace
\end{equation}
with the indices in $\mathcal{V}^{c}=\lbrace v_1,\dots,v_k\rbrace$ in any order. 
\\
\subsubsection{Closure properties}
\begin{figure}[h!]%was scale=0.7
    \centering\fbox{
    \includegraphics[width=0.9\textwidth]{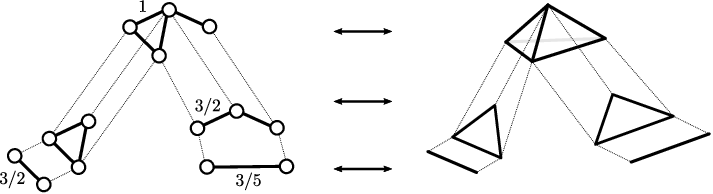}}
    \caption{An example of a graph-Simplex pair and some graphs corresponding to Schur complements of the Laplacian, and faces of the Simplex. By Proposition \ref{Proposition: Gram matrix of face is Schur complement} these sub-objects are again in correspondence. The numbers on the graph indicate the link weights, with unit weights omitted.}
    \label{fig: graph simplex and subsets}
\end{figure}

We now return to the particular case of hyperacute Simplices where the canonical pseudoinverse Gram matrix is Laplacian. In the context of graph Laplacians, the Schur complement 
\begin{equation}\label{eq: definition Schur complement} 
Q/\mathcal{V}^c \triangleq Q_{\mathcal{V}\mathcal{V}} - Q_{\mathcal{V}\mathcal{V}^c}(Q_{\mathcal{V}^c\mathcal{V}^c})^{-1}Q_{\mathcal{V}^c\mathcal{V}}
\end{equation}
is also known as \emph{Kron reduction}, after the foundational work of Gabriel Kron in the study of networks and their reductions \cite{Kron}. For an extensive discussion on Schur complements (Kron reductions) as a tool in electrical circuit and graph theory, we refer to the survey \cite{Dorfler_kron}.
\\
Decomposing the Schur complement of a Laplacian matrix incrementally as in expression \eqref{eq: incremental Schur complement} leads to the following important closure result:
\begin{property}[closure]\label{Property: closure property of Schur complement}
The Schur complement of a Laplacian matrix is again a Laplacian matrix.
\end{property}
\textbf{Proof:} By Property \ref{Property: composition property of Schur complement} and consequently expression \eqref{eq: incremental Schur complement} every Schur complement can be written as a repeated Schur complement with respect to all but a single index $v\in\mathcal{N}$. We will show that for any Laplacian matrix $Q$ and any index $v$, the Schur complement $Q/\lbrace v\rbrace$ is again Laplacian.
We permute the rows and column of $Q$ such that $v$ is in the last position, which gives the block-structure: $Q=\left(\begin{smallmatrix}Q'+\operatorname{diag}(q)&-q\\-q^T&d_v \end{smallmatrix}\right)$, where $Q'$ is the Laplacian of the $(n-1)$-node graph without node $v$, where $(q)_i=w_{vi}$ if $(i,v)\in\mathcal{L}$ and with $d_v$ the degree of $v$. The Schur complement equals $Q/\lbrace v\rbrace=Q'+(\operatorname{diag}(q)-qq^T/d_v)$, where the second term (between brackets) has positive diagonal and non-positive off-diagonal showing that $Q/\lbrace v\rbrace$ has a positive diagonal and non-positive off-diagonal. Furthermore, the matrix is symmetric by construction, and we find that $Q/\lbrace v\rbrace u = Q'u + \operatorname{diag}(q)u -qq^Tu/d_v= 0$ since $Q'u=0$ and $q^Tu=d_v$. Finally, assume the Laplacian $Q/\lbrace v\rbrace$ were reducible, then there exists a bipartition $\mathcal{N}\backslash\lbrace v\rbrace=\mathcal{W}_1\cup\mathcal{W}_2$ such that $(Q/\lbrace v\rbrace)_{ij}=0$ for all $i\in\mathcal{W}_1,j\in\mathcal{W}_2$. However, this would imply that $(Q')_{ij}=0$, i.e. $\mathcal{W}_1$ and $\mathcal{W}_2$ are not connected in $Q'$, as well as $q_iq_j=0$, i.e. the node $v$ can be connected to, say $\mathcal{W}_1$, in $Q$ but not to $\mathcal{W}_2$ at the same time. In other words, this would imply that the partitions $\mathcal{W}_1\cup\lbrace v\rbrace$ and $\mathcal{W}_2$ are not connected in $Q$ which is in contradiction with the fact that $Q$ is Laplacian and thus irreducible. Hence, also $Q/\lbrace v\rbrace$ is irreducible and thus Laplacian.\hfill$\square$
\\
Interestingly, the closure property of Laplacian matrices in combination with Fiedler's identity \eqref{eq: fiedlers identity} yields a surprising closure result for hyperacute Simplices. From Property \ref{Property: face of a Simplex is Simplex} we know that the faces of a hyperacute Simplex are Simplices, but the stronger result that these faces are hyperacute holds as well:
\begin{property}[closure, Fiedler {\cite[Thm. 3.3.2]{Fiedler_book}}]\label{Property: closure property hyperacute Simplices}
The face of a hyperacute Simplex is again a hyperacute Simplex.
\end{property}
\textbf{Proof:} From Proposition \ref{Proposition: Gram matrix of face is Schur complement} and Property \ref{Property: closure property of Schur complement} we have that the canonical pseudoinverse Gram matrix of any face $\mathcal{S}/\mathcal{V}^c$ of a hyperacute Simplex is a Laplacian matrix. Consequently, the face is hyperacute as well by the bijection between Laplacians and hyperacute Simplices, Theorem \ref{Theorem: bijection between Laplacians and Simplices}.\hfill $\square$
\\
Property \ref{Property: closure property hyperacute Simplices} is surprising since it implies that the $n(n-1)/2$ angle constraints between the dihedral angles of a Simplex somehow also influence the dihedral angles between \emph{all} faces of this Simplex. The bijection Theorem \ref{Theorem: bijection between Laplacians and Simplices} combined with the fact that the face relation of a Simplex translates to a Schur complement relation, which maps Laplacians to Laplacians, shows how exactly this (top-level) Simplex constraint on the angles is inherited by the (lower-level) faces. The interrelated `nested structures' of graphs and Simplices is illustrated by an example in Figure \ref{fig: graph simplex and subsets}.
%%%
%
%%%
\section{A geometric proof of the resistance distance}
\label{Section: geometric proof of resistance distance}
Our discussion of the relation between Laplacians and hyperacute simplices now enables a geometric and intuitive proof of the fact that the effective resistance is a metric function; we recall that a function $d:\mathcal{N}\times\mathcal{N}\rightarrow \mathbb{R}_+$ is \emph{metric} if it satisfies the following criteria:
\begin{align*}
    &\text{(i)$_m$~}d(x,y)=0\text{~if and only if~}x=y\\
    &\text{(ii)$_m$~}d(x,y)=d(y,x)\text{~for all~}x,y\\
    &\text{(iii)$_m$~}d(x,y)+d(y,z)\geq d(x,z)\text{~for all~}x,y,z
\end{align*}
for any three elements $x,y,z\in\mathcal{N}$ in the set. Properties (i)$_m$-(ii)$_m$ are usually easily confirmed, while property (iii)$_m$, the \emph{triangle inequality}, is typically harder to verify for a candidate metric function $d$. 
\\
In our case, we are interested to show that the effective resistance $d(i,j)=\omega_{ij}$ is a distance function. To start, we have already shown the following related result:
\begin{lemma}\label{Lemma: square root of resistance is Euclidean metric of hyperacute Simplex}
The square root of the effective resistance $\sqrt{\omega}:(i,j)\mapsto\sqrt{\omega_{ij}}$ is a metric function. Moreover, $\sqrt{\omega}$ is a Euclidean metric which embeds the elements in $\mathcal{N}$ as the vertices of a hyperacute Simplex.
\end{lemma}
\textbf{Proof:} This Lemma follows from expression \eqref{eq: effective resistance as edge lengths} for the effective resistance and Lemma \ref{Lemma: Laplacians are Gram matrices}.\hfill$\square$
\\
From the fact that $\sqrt{\omega}$ is metric, it follows by the square relation that the effective resistance $\omega$ satisfies properties (i)$_m$-(ii)$_m$ as well. It thus remains to verify the triangle inequality on all triples $(\omega_{ij},\omega_{jk},\omega_{ik})$. By Lemma \ref{Lemma: square root of resistance is Euclidean metric of hyperacute Simplex}, every such triple corresponds to the squared edge-lengths of a triangular face of the hyperacute Simplex corresponding to $\sqrt{\omega}$. By the closure property of hyperacute simplices, every such triangular face is a hyperacute triangle and we have that \emph{every triple $(\omega_{ij},\omega_{jk},\omega_{ik})$ corresponds to the squared edge-lengths of a hyperacute triangle}. Finally, by basic trigonometry (for instance, the cosine rule) we have that a triangle with squared edge-lengths $(\omega_{ij},\omega_{jk},\omega_{ik})$ and a non-obtuse angle between the edges with lengths $\omega_{ij}$ and $\omega_{jk}$ satisfies
$$
\omega_{ij}+\omega_{jk}\geq \omega_{ik}.
$$
Since this holds for any of the three angles, and for any triple of effective resistances, the triangle inequality is satisfied in general and we have:
\begin{theorem}[\cite{Gvishiani_Gurvich} and \cite{Klein}]\label{Theorem: resistance is distance}
The effective resistance $\omega:(i,j)\mapsto\omega_{ij}$ is a metric function.
\end{theorem}
\textbf{Proof:} We summarize the derivation above; by Lemma \ref{Lemma: square root of resistance is Euclidean metric of hyperacute Simplex} the effective resistances are squared edge lengths of a hyperacute Simplex. By the closure property of hyperacute Simplices, any triangular face is hyperacute as well such that its squared edge lengths, which are the effective resistances, satisfy the triangle inequality. \hfill $\square$
\\
The geometric proof described above highlights the fact that the effective resistance being metric is actually just a manifestation of a much richer structure. In a certain sense, Lemma \ref{Lemma: square root of resistance is Euclidean metric of hyperacute Simplex} which states that the (square root) effective resistance determines a hyperacute Simplex better captures this structure. One could in fact write down the `non-obtuse dihedral angle' constraint for any $v$-dimensional face in terms of the effective resistances, which would yield a set of (non-trivial) inequalities that the effective resistances must satisfy; the triangle inequality is then just one example of such inequalities, obtained from the non-obtuseness of triangular faces. This perspective is also explored by Klein in \cite{Klein_wiener}. Another important remark is that it is not sufficient for a simplex to just have hyperacute triangular faces in order for its squared edge lengths to determine an `effective resistance metric', but that the hyperacute inequalities need to be satisfied for all faces simultaneously; the following example illustrates this fact.
\\
\textit{Example:} below we give an example of a non-hyperacute Simplex $\mathcal{S}$ with hyperacute triangular facets, whose squared edge-lengths thus determine a metric. The simplex $\mathcal{S}$ and its faces have the following canonical pseudoinverse Gram matrices:
\begin{align*}\label{eq: counterexample}
&M^\dagger(\mathcal{S}) = 18\begin{pmatrix}9&1&-5&-5\\1&9&-5&-5\\-5&-5&15&-5\\-5&-5&-5&15 \end{pmatrix} \text{~with~} \\
&M^\dagger(S/\lbrace 1\rbrace)=M^\dagger(S/\lbrace 2\rbrace)=
20\begin{pmatrix} 8&-4&-4\\-4& 11 &-7\\-4&-7& 11\end{pmatrix}\\
&M^\dagger(S/\lbrace 3\rbrace)=M^\dagger(S/\lbrace 4\rbrace)=
12\begin{pmatrix}11&-1&-10\\-1&11 &-10\\-10&-10&20 \end{pmatrix}
\end{align*}
Since $(M^\dagger(\mathcal{S}))_{12}>0$ the angle between $S_{\lbrace 1\rbrace^{c}}$ and $S_{\lbrace 2\rbrace^{c}}$ is obtuse and hence the Simplex $\mathcal{S}$ is not hyperacute. However, all triangular faces have pseudoinverse Gram matrices with $[M^\dagger(\mathcal{S}/\lbrace i\rbrace)]_{ab}< 0$ for $a\neq b$ which are Laplacian and thus correspond to hyperacute triangles. 
\\~\\
\textbf{Conclusion:} The example above is a good reflection of the key message of this article: \emph{the effective resistance is more than just a distance}; it reflects the geometric structure of graphs as a hyperacute simplex whose properties translate to properties of the Laplacian via Fiedler's algebraic identity. Moreover, the key role of the Schur complement, which maps Laplacians to Laplacians and the corresponding hyperacute simplices to hyperacute simplices is clearly highlighted in our setup for the geometric proof of the resistance distance.
\\~\\
\textbf{Outlook:} While not further developed in this document, our description of the relation between graphs via the Schur complement and Simplices via the face relation, which are both closed and composable, provides the necessary setup to define a \emph{category of graphs} $\mathbf{G}$ and a \emph{category of hyperacute} Simplices $\mathbf{S}$ (see e.g. \cite{Leinster_cat_th_book} for an introduction to category theory). The objects in $\mathbf{G}$ are finite connected graphs with finite, non-negative link weights and the morphisms between graphs follow from the Schur complement. The objects in $\mathbf{S}$ are finite, non-degenerate hyperacute Simplices and the morphisms between Simplices correspond to the face relation. In this setup, the bijection of Theorem \ref{Theorem: bijection between Laplacians and Simplices} between the objects of both categories is in fact a \emph{functorial} relation between categories $\mathbf{G}$ and $\mathbf{S}$ due to Proposition \ref{Proposition: Gram matrix of face is Schur complement} which implies that the diagram below commutes:
\[ \begin{tikzcd}
G \arrow[leftrightarrow]{r}{f} \arrow[swap]{d}{\phi} & \mathcal{S} \arrow{d}{\psi} \\%
G' \arrow[leftrightarrow]{r}{g}& \mathcal{S}'
\end{tikzcd}\text{~with~$G,G'\in\mathbf{G}$ and $\mathcal{S},\mathcal{S}'\in\mathbf{S}$.}
\]
In other words, there is not only a bijection between graphs and Simplices, but this bijection also respects the interconnection structure (i.e. morphisms) between the respective objects themselves. Consequently, a stronger version of Theorem \ref{Theorem: bijection between Laplacians and Simplices} would say that \emph{the categories $\mathbf{G}$ and $\mathbf{S}$ are equivalent}, which is a more complete characterization of the relation between graphs and simplices discovered by Fiedler. As an outlook, we might hope that this abstract categorical perspective provides a new stepping stone for further developments in the theory of graphs, Laplacians, effective resistances and hyperacute simplices. In particular, we believe that the above described structure `sits inside' the larger categorical framework for passive linear circuits described in \cite{Baez_Fong} by Baez and Fong and that a specialization of their results to our setting and, similarly, a generalization of our results to their broader setting would be an interesting way forward.

%%%%
%%%%
%%%%
%%%%%%%%%%%%%%%%%%
% ACKNOWLEDGEMENTS
%%%%%%%%%%%%%%%%%%
\section*{Acknowledgements}
This work was initiated at TUDelft under supervision of Piet Van Mieghem. While finishing and writing up, the author was supported by The Alan Turing Institute under the EPSRC grant EP/N510129/1. KD is grateful to Piet Van Mieghem, Renaud Lambiotte, Marc Homs-Dones and an anonymous referee for their helpful comments and suggestions and to KVL for their support.
\bibliographystyle{elsarticle-num} 
\bibliography{bibliography.bib}

\end{document}